\documentclass[12pt]{amsart}

\usepackage{amssymb,amsfonts,amsthm,amsmath,amscd}
\usepackage{graphics}
\usepackage{graphicx}
\usepackage{supertabular}
\usepackage{tikz}
\usepackage{multirow}
\usepackage{enumerate}

\usetikzlibrary{calc, shapes, backgrounds }
\usepackage{pgfplots,pgfplotstable}
\pgfplotsset{width=7cm}

\newtheorem{theorem}{Theorem}
\newtheorem{lemma}{Lemma}

\newtheorem{corollary}{Corollary}

\theoremstyle{definition}

\newcommand{\ga}{\alpha}
\newcommand{\gb}{\beta}
\newcommand{\gc}{\gamma}

\def\G{\mathcal{G}}
\def\T{\mathcal{T}}
\def\P{\mathcal{P}}
\def\N{\mathcal{N}}

\def\Z{\mathcal{Z}}
\def\ZZP{\mathbb{Z}_{\geq}}
\def\ZZ{\mathbb{Z}}

\def\cH{\mathcal{H}}

\begin{document}

\title{Tetris Hypergraphs and Combinations of Impartial Games}

\author{Endre Boros}
\address{MSIS and RUTCOR, RBS, Rutgers University,
100 Rockafeller Road, Piscataway, NJ 08854}
\email{endre.boros@rutgers.edu}

\author{Vladimir Gurvich}
\address{MSIS and RUTCOR, RBS, Rutgers University,
100 Rockafeller Road, Piscataway, NJ 08854; \\
Dep. of Computer Sciences, National Research University,
Higher School of Economics (HSE), Moscow}
\email{vladimir.gurvich@rutgers.edu}

\author{Nhan Bao Ho}
\address{Department of Mathematics and Statistics, La Trobe University, Melbourne, Australia 3086}
\email{nhan.ho@latrobe.edu.au, nhanbaoho@gmail.com}

\author{Kazuhisa Makino}
\address{Research Institute for Mathematical Sciences (RIMS)
Kyoto University, Kyoto 606-8502, Japan}
\email{e-mail:makino@kurims.kyoto-u.ac.jp}

\author{Peter Mursic}
\address{MSIS and RUTCOR, RBS, Rutgers University,
100 Rockafeller Road, Piscataway, NJ 08854}
\email{peter.mursic@rutgers.edu}


\thanks{The first author gratefully acknowledges the support by Kyoto University during his visit in the Fall of 2015.
The second author was partially funded by
the Russian Academic Excellence Project '5-100'.
The fourth author was partially supported by KAKENHI Grant Numbers 24106002 and 26280001.}

\begin{abstract}
The Sprague-Grundy (SG) theory reduces the sum of
impartial games to the classical game of $NIM$.
We generalize the concept of sum and introduce
$\cH$-combinations of impartial games for any hypergraph $\cH$.
In particular, we introduce the game $NIM_\cH$ which is the $\cH$-combination of single pile $NIM$ games.
An impartial game is called SG decreasing if its SG value is
decreased by every move.
Extending the SG theory, we reduce the $\cH$-combination of SG decreasing games to $NIM_\cH$.
We call $\cH$ a Tetris hypergraph if $NIM_\cH$ is SG decreasing.
We provide some necessary and some sufficient conditions for a hypergraph to be Tetris.
\end{abstract}

\maketitle

\section{Introduction}
\label{s1}

An impartial game can be modeled by a directed graph
$\Gamma = (X, E)$, in which a vertex
$x \in X$  represents a {\em position}, while a directed edge
$(x, x') \in E$  represents a {\em move} from position $x$  to  $x'$,
which we will also denote by $x \to x'$.

The graph  $\Gamma$  may be infinite, but we will always
assume that any sequence of successive moves
(called a {\em play})  $x \to x'$, $x' \to x''$, $\cdots$  is  finite.
In particular, this implies that $\Gamma$  has no directed cycles.
The game is played by two players with a token placed at an initial position.
They alternate in moving the token along the directed edges of the graph.
The game ends when the token reaches a {\em terminal}, that
is, a vertex with no outgoing edges.
The player who made the last move wins, equivalently, the one who is out of moves, loses.

In this paper we consider only impartial games and call them simply games.

The basic concepts and definitions related to impartial games
will be  briefly summarized in the next section.
We refer the reader to \cite{Alb07, BCG01-04} for more details.

It is known that the set of positions of a game $\Gamma$ can
uniquely be partitioned into sets of winning and loosing positions.
Every move from a loosing position goes to a winning one, while
from a winning position we always have a move to a loosing one; see Section \ref{s2}.
This partition shows how to win the game, whenever possible.
The so-called Sprague-Grundy (SG) function $\G_\Gamma:X\to \ZZ_\geq$ is
a generalization of the above partition; see Section \ref{s2}.
Namely, $\G_\Gamma(x)=0$ if and only if $x$ is a loosing position.

Given $n$ games $\Gamma_i$, $i=1,...,n,$
their {\em sum} $\Gamma = \Gamma_1 + \cdots +\Gamma_n$
is the game in which players on their turn choose one of the games and make a move in it.
To play optimally the sum, it is not enough to know the winning-loosing partitions
of all the games $\Gamma_i$, $i=1,...,n$.
Sprague and Grundy \cite{Spr35, Spr37, Gru39} resolved this problem.
They proved that the SG function of the sum can easily be computed from the SG functions of the summands
and thus we can compute the winning-loosing partition of the sum; see Section \ref{s2}.

A classical example for sum is $NIM_n$.
Given $n$ piles of stones, a move consists of
choosing a nonempty pile and taking some positive number of stones from it.
By this definition, $NIM_n$ is the sum of $n$ single pile NIM games.
Bouton \cite{Bou901} described the winning-loosing partition of this game.

In this paper we propose a generalization of the notion of sum of games.
Given $n$ games $\Gamma_i$, $i\in V=\{1,...,n\}$, as above, and
a hypergraph $\cH \subseteq 2^V$, we introduce the $\cH$-{\em combination} $\Gamma_\cH$ of the given games.
In $\Gamma_\cH$ players on their turn choose one of the hyperedges
$H \in \cH$ and make a move in all $\Gamma_i$, $i\in H$.
As before, the player who cannot do this is the looser.

The sum of $n$ games is an $\cH$-combination with $\cH=\{\{1\},...,\{n\}\}$.
If $\cH=\{V\}$, we call the $\cH$-combination the \emph{product} of the given games.
Let us note that these and a few similar operations were considered under different names in the literature.
For instance, Smith \cite{Smi66} calls the sum a disjunctive compound, the product a conjunctive compound;
he also mentions selective compounds, which are $\cH$-combinations with $\cH=\{S\subseteq V \mid S \neq \emptyset\}$.

Moore \cite{Moo910} introduced a generalization of $NIM_n$, which we will denote $NIM_{n,k}^\leq$.
Given $n$  piles of stones, a move in $NIM_{n, k}^\leq$ consists of choosing $\ell$ nonempty piles,
$1 \leq \ell \leq k,$ and taking a positive number of stones from each of the chosen piles.
Moore \cite{Moo910} described the winning-loosing partition of $NIM_{n,k}^\leq$, in other words, positions of SG value $0$. Jenkyns and Mayberry \cite{JM80} described the set of positions of $NIM_{n,k}^\leq$ in which the SG value is $1$ and got an explicit formula for the SG function in case  $k = n-1$.
The game $NIM_{n,k}^\leq$ can also be seen as an $\cH$-combination
with $\cH = \cH_{n,k}^\leq = \{S\subseteq V \mid 1 \leq |S| \leq k\}$.

In \cite{BGHMM15} another generalization $NIM_{n,k}^=$  of $NIM_n$  was introduced, which
is also an $\cH$-combination with  $\cH=\cH_{n,k}^= = \{S\subseteq V \mid |S|=k\}$, and
an explicit formula was given for the SG function when $2k \geq n$.

Given a hypergraph $\cH \subseteq 2^V$, let us define $NIM_{\cH}$ as
the $\cH$-combination of $n$ single pile $NIM_1$ games. Note that the family of $NIM_\cH$ games is closed under hypergraph combinations. The above cited generalizations of $NIM$ all belong to the family of $NIM_\cH$ games.

In all of the above cases, when the SG functions are known, we have the following equality:

\begin{equation}\label{e1}
\G_{\Gamma_\cH} ~=~ \G_{NIM_\cH} \left(\G_{\Gamma_1},...,\G_{\Gamma_n}\right).
\end{equation}

In this paper we would like to find other classes of games for which equality \eqref{e1} holds.

A game $\Gamma$ is called \emph{SG decreasing} if the SG value is strictly decreased by every move. The single pile $NIM$ is the simplest example of an SG decreasing game. Another example is $NIM_{n,k}^=$ when $n<2k$, \cite{BGHMM15}.
We call a hypergraph $\cH$ \emph{Tetris} if $NIM_\cH$ is an SG decreasing game.

Given a game $\Gamma$ and a position $x$ of it,
we denote by $\T_\Gamma(x)$ the length of the longest play starting at $x$ and
call $\T_\Gamma$ the \emph{Tetris} function of $\Gamma$.
It is easy to show that $\Gamma$ is SG decreasing if and only if $\G_\Gamma = \T_\Gamma$; see Section \ref{s4}
This implies that from every nonterminal position of an SG decreasing game one can win by a single move.

\begin{theorem}
\label{t1}
Equality \eqref{e1} holds for arbitrary hypergraph
$\cH \subseteq 2^V$ and SG decreasing games $\Gamma_i$, $i\in V$.
\end{theorem}

Let us remark that equality \eqref{e1} does not always hold.
Consider, for example, $\Gamma_1=NIM_1$, $\Gamma_2=NIM_2$ and  $\cH=\{\{1,2\}\}$.
In this case $\Gamma_2$ is not SG decreasing and equality \eqref{e1} may fail.

While computing the SG function for games seems to be very hard, in general,
the above theorem allows us to outline new cases when the problem is tractable.
In \cite{BGHMM16} we introduce a special family of hypergraphs
for which $\G_{NIM_\cH}$ can be described by a closed formula.
This result combined with Theorem \ref{t1} provides new families of games for which
the SG function can be expressed by an explicit formula
via the Tetris functions $\T_{\Gamma_i}$, $i \in V$.
For instance, the above introduced hypergraphs $\cH_{2k,k}^=$ for $k \geq 2$ and
$\cH_{k+1,k}^\leq$ for $k \geq 2$  appear to be such families.

Let us add that sometimes even for very small games we do not know how to compute the SG function; for instance, for the games $NIM_\cH$ with
$\cH=\{\{1,2\},\{1,3\},\{1\},\{2\},\{3\}\}$ \cite{BGHM15}, or
$NIM_{4,2}^\leq$ \cite{JM80, BGHM15}, or $NIM_{5,2}^=$ \cite{BGHMM15}.
In the latter case we cannot even describe the winning-loosing partition.

\medskip

If we replace in equality \eqref{e1} the SG function by the Tetris function we get always equality.

\begin{theorem}
\label{t2}
Given games $\Gamma_i$, $i\in V$, and a hypergraph $\cH\subseteq 2^V$ we have the equality
\begin{equation}\label{e2}
\T_{\Gamma_\cH} ~=~ \T_{NIM_\cH} \left(\T_{\Gamma_1},...,\T_{\Gamma_n}\right).
\end{equation}
\end{theorem}

The above two theorems immediately imply the following statement.

\begin{corollary}
\label{c-superposition}
If $\cH$ is a Tetris hypergraph then
 $\cH$-combination of SG decreasing games is SG decreasing.
 In particular, a Tetris combination of Tetris hypergraphs is Tetris.
\end{corollary}

\bigskip

While recognizing if a given hypergraph is Tetris is a hard decision problem,
we can provide a necessary and sufficient condition for hypergraphs of dimension at most $3$,
where the dimension $dim(\cH)$ of $\cH$ is the size of the largest hyperedge in $\cH$.
For a hypergraph $\cH \subseteq 2^V$ and a subset $S\subseteq V$ we denote by $\cH_S$ the subhypergraph of $\cH$ induced by $S$
\[
\cH_S ~=~ \{H\in \cH\mid H\subseteq S\}.
\]

\begin{theorem}
\label{t3}
A hypergraph $\cH$ of dimension at most $3$ is Tetris if and only if
\begin{equation}\label{intersection}
\forall S\subseteq V ~\text{ with }~ \cH_S\neq\emptyset ~~\exists H\in\cH_S ~\text{ such that }~ H\cap H'\neq\emptyset ~\forall H'\in\cH_S.
\end{equation}
\end{theorem}

Let us remark that computing the Tetris function value for $NIM_\cH$ is an NP-hard problem,
even for hypergraphs of dimension at most $3$.
Furthermore, computing the SG function value is also NP-hard for this family of games.

\section{Definitions and notation}
\label{s2}

It is not difficult to characterize
the winning strategies in a game $\Gamma = (X, E)$.
The subset  $\P \subseteq X$  of the loosing positions
is uniquely defined by the following two properties:

\begin{itemize}
\item[(IND)]
$\P$ is {\em independent}, that is,
for any  $x \in \P$  and move $x \to x'$  we have  $x' \not\in \P$;
\item[(ABS)]
$\P$ is {\em absorbing}, that is, for any  $x \not\in \P$
there is a move  $x \to x'$  such that $x' \in P$.
\end{itemize}

It is easily seen that the
set $\P$ can be obtained by the following simple recursive algorithm \cite{NM44}:
include in  $\P$  each terminal of  $\Gamma$;
include in  $X \setminus \P$  every position  $x$  of  $\Gamma$
from which there is a move  $x \to x'$  to a terminal  $x'$;
delete from $\Gamma$  all considered positions and repeat.

It is also clear that
any move  $x \to x'$  of a player to a $\P$-{\em position}
$x' \in \P$  is a winning move.
Indeed, by (IND), the opponent must leave  $\P$
by the next move, and then, by (ABS), the player can reenter  $\P$.
Since, by definition, all plays of  $\Gamma$ are finite
and, by construction, all terminals are in $P$, sooner or later
the opponent will be out of  moves.

In combinatorial game theory positions  $x \in \P$  and  $x \not\in \P$
are usually called a $\P$- and $\N$-positions, respectively.
The next player wins in an  $\N$-position, while
the previous one wins in a $\P$-position.

\smallskip



By definition, $NIM_n$ is the sum of $n$ games, each of which
(a single pile $NIM_1$) is trivial.
Yet, $NIM_n$ itself is not. It was solved by Bouton in his seminal paper \cite{Bou901} as follows.
The $NIM$-sum  $x_1 \oplus \cdots \oplus x_n$ of nonnegative integers
is defined as the bitwise binary sum. For example,

\smallskip

$3 \oplus 5 = 011_2 \oplus 101_2 = 110_2 = 6$, $3 \oplus 6 = 5$, $5 \oplus 6 = 3$, and $3 \oplus 5 \oplus 6 = 0$.

\smallskip
\noindent
It was shown in \cite{Bou901} that
$x = (x_1, \ldots, x_n)\in \ZZP^n$ is a $\P$-position of $NIM_n$
if and only if  $x_1 \oplus \cdots \oplus x_n = 0$.

\smallskip

To play the sum  $\Gamma = \Gamma_1 + \Gamma_2$, it is not sufficient
to know  $\P$-positions of  $\Gamma_1$ and $\Gamma_2$, since
$x = (x^1, x^2)$  may be a  $\P$-position of  $\Gamma$  even when
$x^1$ is not a $\P$-position of $\Gamma_1$  and
$x^2$ is not a $\P$-position of $\Gamma_2$.
For example,  $x = (x_1, x_2)$  is a $\P$-position of the two pile
$NIM_2$ if and only if  $x_1 = x_2$, while
only  $x_1 =0$  and  $x_2 =0$   are
the unique $\P$-positions of the corresponding single pile games.
To play the sums  we need the concept of the Sprague-Grundy
(SG) function, which
is a refinement of the concept of $\P$-positions.

Given a finite subset 
$S \subseteq \ZZP$, let
$mex(S)$  (the \emph{minimum excluded value}) be
the smallest  $k \in \ZZP$  that is not in $S$.
In particular, $mex(\emptyset)=0$, by the definition.

Given an impartial game $\Gamma = (X, E)$, the SG function
$\G_\Gamma : X \to \ZZP$  is defined recursively, as follows:
$\G_\Gamma(x) = 0$  for any terminal  $x$ and, in general,
$\G_\Gamma(x) = mex (\{\G_\Gamma(x') \mid x \to x'\})$.

It can be seen easily that the following two properties define the SG function uniquely.

\begin{itemize} \itemsep0em
\item[(1)] No  move keeps the SG value, that is, $\G_\Gamma(x) \neq \G_\Gamma(x')$  for any move $x \to x'$.
\item[(2)]  The SG value can be arbitrarily (but strictly) reduced by a move, that
is, for any integer  $v$  such that  $0 \leq v < \G(x)$
there is a move  $x \to x'$  such that  $\G_\Gamma(x') = v$.
\end{itemize}

The definition of the SG function implies several other important properties:

\begin{itemize} \itemsep0em
\item[(3)] The $\P$-positions are exactly the zeros of the SG function:
$\G_\Gamma(x) = 0$  if and only if  $x$  is a $\P$-position of $\Gamma$.
\item[(4)] The SG function of $NIM_n$ is the $NIM$-sum of
the cardinalities of its piles, that is,
$\G_{NIM_n}(x) = x_1  \oplus \cdots \oplus x_n$  for all $x = (x_1, \ldots, x_n)\in \ZZP^n$; see \cite{Bou901,Spr35,Spr37,Gru39}.
\item[(5)] In general, the SG function of the sum of  $n$  games is
the $NIM$-sum of the  $n$  SG functions of the summands.
More precisely, let
$\Gamma = \Gamma_1 + \ldots + \Gamma_n$  be the sum of $n$ games and
$x = (x^1, \ldots, x^n)$  be a position of  $\Gamma$, where
$x^i$  is a position of
$\Gamma_i$ for $i \in V=\{1,...,n\}$, then
$\G_\Gamma(x) = \G_{\Gamma_1}(x^1) \oplus \cdots \oplus \G_{\Gamma_n}(x^n)$; see \cite{Spr35,Spr37,Gru39}.
\end{itemize}


SG theory shows that playing a sum of games
$\Gamma = \Gamma_1 + \cdots + \Gamma_n$
may be effectively replaced by $NIM_n$
in which each summand game $\Gamma_i$ is replaced
by a pile of  $x_i = \G_{\Gamma_i}(x)$  stones, for $i\in V$.

\section{Hypergraph Combinations of Games}
\label{s3}

Given games $\Gamma_i=(X_i,E_i)$, $i\in V=\{1,...,n\}$, and a hypergraph $\cH\subseteq 2^V$, we define the $\cH$-combination $\Gamma_\cH=(X,E)$ of these games by setting
\[
X ~=~ \prod_{i\in V} X_i, \text{ and }
\]
\[
E ~=~ \left\{(x,x')\in X\times X ~ \left| ~\exists ~ H\in\cH \text{ such that }
\begin{array}{ll}
(x_i,x_i')\in E_i & \forall i\in H, \\
x_i=x_i' & \forall i\not\in H
\end{array}
\right.\right\} .
\]

Since the combination game $NIM_\cH$ plays a special role in our statements, we introduce a simplified notation for the rest of the paper. Namely, we denote by $\G_\cH$ the SG function, by $\T_\cH$ the Tetris function, and by $\P_\cH$ the set of $P$-positions of $NIM_\cH$.

\subsection*{Proof of Theorem \ref{t1}}

Consider the $\cH$-combination $\Gamma_\cH=(X,E)$ as defined above, and show that the function defined by \eqref{e1} satisfies the defining properties of the SG function.

First, consider a position $x=(x_1,...,x_n)\in X$.  Let us denote by $g=g(x)=(\G_{\Gamma_i}(x_i)\mid i\in V)\in \ZZP^V$ the vector of SG values in the $n$ given games. Note that $g\in \ZZP^V$ is a position in the game $NIM_\cH$.

Let us denote by
\[
f(x_1,...,x_n)=\G_\cH(g(x))
\]
the function defined by the right hand side of \eqref{e1}.

Consider first a move $(x,x')\in E$, where
$(x_i,x_i')\in E_i$ for $i\in H$ for some hyperedge $H\in \cH$. By the definition of $E$ we must have $x_i'=x_i$ for all $i\not\in H$. Denote by $g'\in \ZZP^V$ the corresponding vector of SG values. Note that $g'_i<g_i$ for $i\in H$ since $\Gamma_i$ is an SG decreasing game for all $i\in V$, and $g_i'=g_i$ for all $i\not\in H$ since $x_i=x_i'$ for these indices. Consequently, $g\to g'$ is a move in $NIM_\cH$, and therefore $f(x)=\G_\cH(g)\neq \G_\cH(g')=f(x')$. Thus we proved that every move in $\Gamma_\cH$ changes the value of function $f$.

Next, let us consider an integer $0\leq v < f(x)$. We are going to show that there exists a move $x\to x'$ in $\Gamma_\cH$ such that $f(x')=v$. Let us consider again the corresponding integer vector $g=g(x)\in \ZZP^V$, for which we have $f(x)=\G_\cH(g)$. By the definition of the SG function of $NIM_\cH$, there must exists a move $g\to g'$ such that $\G_\cH(g')=v$. Assume that this move is an $H$-move for some $H\in \cH$, that is that $g'_i<g_i$ for $i\in H$ and $g'_i=g_i$ for $i\not\in H$. Then we have $\G_{\Gamma_i}(x_i)=g_i > g_i'$ for all $i\in H$, and thus we must have moves $x_i\to x_i'$ in $\Gamma_i$, $i\in H$ such that $\G_{\Gamma_i}(x_i')=g_i'$ for all $i\in H$. Then with $x_i'=x_i$ for $i\not\in H$, we get that $x\to x'$ is a move in the $\cH$-combination, and $f(x')=v$. Thus we proved that each smaller SG value can be realized by a move in the combination game.

The above arguments can be completed by an easy induction to show that $\G_{\Gamma_\cH}(x)=f(x)$ for all $x\in X$.
\qed

\bigskip

Since the Tetris function is the length of a longest path in the directed graph of the game, it is uniquely defined by the following three properties:
\begin{itemize} \itemsep0em
\item[(a)] Every move decreases its value.
\item[(b)] If it is positive in a position, then there exists a move from that position that decreases it by exactly one.
\item[(c)] It takes value zero at every terminal.
\end{itemize}

\smallskip

\subsection*{Proof of Theorem \ref{t2}}

Similarly to the proof of Theorem \ref{t1}, we shall show that the function defined by the right hand side of \eqref{e2} satisfies properties (a), (b) and (c) above.

Consider a position $x=(x_1,...,x_n)\in X$, and denote by $t=t(x)=(\T_{\Gamma_i}(x_i)\mid i\in V)\in \ZZP^V$ the vector of Tetris values in the $n$ given games. Note that $t$ is a position in the game $NIM_\cH$.

Let us denote by
\[
h(x_1,...,x_n)=\T_\cH(t(x))
\]
the function defined by the right hand side of \eqref{e2}.

Consider first a move $(x,x')\in E$, where
$(x_i,x_i')\in E_i$ for $i\in H$ for some hyperedge $H\in \cH$. By the definition of $E$ we must have $x_i'=x_i$ for all $i\not\in H$. Denote by $t'\in \ZZP^V$ the corresponding vector of Tetris values, and note that $t'_i<t_i$ for $i\in H$ since $T_{\Gamma_i}$ satisfies property (a) for all $i\in V$, and $t_i'=t_i$ for all $i\not\in H$ since $x_i=x_i'$ for these indices. Consequently, $t\to t'=t(x')$ is a move in $NIM_\cH$, and therefore $h(x)=\T_\cH(t(x))> \T_\cH(t(x'))=h(x')$, since $\T_\cH$ satisfies property (a).
Thus we proved that every move in $\Gamma_\cH$ decreases the value of function $h$.

Consider next an arbitrary position $x\in X$ such that $0 < h(x)=\T_\cH(t(x))$. Since $\T_\cH$ satisfies property (b), there exists a move $t(x)\to t'$ in $NIM_\cH$ such that $\T_\cH(t')=\T_\cH(t(x))-1$. Then, by the definition of $NIM_\cH$ we must have $H=\{i\in V\mid t_i>t_i'\}\in \cH$. Since $\T_{\Gamma_i}$
satisfies property (b), there must exist moves $x_i\to x_i'$ such that $\T_{\Gamma_i}(x_i')=\T_{\Gamma_i}(x_i)-1=t_i-1$ for $i\in H$. Define $x_i'=x_i$ for $i\not\in H$. Then we have
$\T_\cH(t(x))-1\geq \T_\cH(t(x'))\geq \T_\cH(t')$ by the definition of $NIM_\cH$. Consequently we have $h(x')=h(x)-1$.

Finally, to see property (c), let us consider a terminal position $x\in X$ and its corresponding Tetris value vector $t(x)$. By the definition of $NIM_\cH$ this is a terminal position if and only if $\{i\in V\mid t_i=0\}$ intersects all hyperedges of $\cH$, in which case we must have $h(x)=T_\cH(t(x))=0$.
\qed

\section{Tetris Hypergraphs}
\label{s4}

\subsection{A necessary condition for $\T_\cH=\G_\cH$}\label{s-necessary-dim}

Let us start by observing that for every game $\Gamma=(X,E)$ and position $x\in X$ we have the inequality
\begin{equation}\label{trivi-UB}
\G_\Gamma (x) ~\leq~ \T_\Gamma (x).
\end{equation}

Let us continue with some basic properties of Tetris hypergraphs.

\begin{lemma}\label{l2-gap}
Given a hypergraph $\cH\subseteq 2^V$, the following three statements are equivalent:
\begin{itemize} \itemsep0em
\item[\rm{(i)}] $\cH$ is a Tetris hypergraph;
\item[\rm{(ii)}]  $\G_{\cH}=\T_{\cH}$;
\item[\rm{(iii)}] for all positions $x\in\ZZP^V$ and for all integers $0\leq v<\T_{\cH}(x)$ we have a move $x\to x'$ in $NIM_\cH$ such that $\T_{\cH}(x')=v$.
\end{itemize}
\end{lemma}

\proof
These equivalences follow directly from the definitions of Tetris hypergraphs, Tetris and SG functions, and SG decreasing games.
\qed

For a hyperedge $H\in\cH$, and a position $x\in\ZZP^V$, we call a move $x\to x'$ in $NIM_\cH$ an $H$-\emph{move} if $\{i\in V\mid x_i'<x_i\}=H$.

For a subset $H\subseteq V$ we denote by $\chi(H)$ its characteristic vector.
For positions $x\geq \chi(H)$ we shall consider two special $H$-moves from $x$:
\begin{description}
\item[Slow $H$-move:] $x\to x^{s(H)}$ defined by $x^{s(H)}_i=x_i-1$ for $i\in H$, and $x^{s(H)}_i=x_i$ for $i\not\in H$, that is by decreasing every coordinate in $H$ by exactly one.
\item[Fast $H$-move:] $x\to x^{f(H)}$ defined by $x^{f(H)}_i=0$ for $i\in H$, and $x^{f(H)}_i=x_i$ for $i\not\in H$, that is by decreasing the size of every coordinate in $H$ to zero.
\end{description}

Let us associate to a hypergraph $\cH\subseteq 2^V$ the set of positions $\Z_\cH\subseteq \ZZP^V$ which have zero Tetris value:
\[
\Z_\cH ~=~ \{x\in \ZZP^V\mid T_\cH(x)=0\}.
\]
Obviously, we have
\begin{equation}\label{trivi-UB-2}
\Z_\cH ~\subseteq~ \P_\cH,
\end{equation}
since there is no move from $x$ by the definition of the Tetris function.
We shall show next that in fact all P-positions of $NIM_\cH$ are in $\Z_\cH$
if and only if condition \eqref{intersection} holds.

\begin{theorem}
\label{t-P-Z-equivalence}
Given a hypergraph $\cH\subseteq 2^V$, $\emptyset\not\in\cH$, we have $\Z_\cH=\P_\cH$  if and only if $\cH$ satisfies property \eqref{intersection}.
\end{theorem}

\proof
By \eqref{trivi-UB-2} we always have $\Z_{\cH}\subseteq \P_{\cH}$.

Assume first that we also have $\P_{\cH}\subseteq \Z_{\cH}$, and consider a subset $S\subseteq V$ such that $\cH_S\neq\emptyset$. For a position $x\in \ZZP^V$ we denote by $supp(x)=\{i\mid x_i>0\}$ the set of its support. Let us then choose a position $x\in\ZZP^V$ such that $supp(x)=S$. Since $\cH_S\neq\emptyset$, we have $\T_{\cH}(x)>0$ implying $\G_{\cH}(x)>0$ by our assumption. Then, by the definition of the SG function we must have a hyperedge $H'\in \cH$ and an $H'$-move $x\to x'$ such that $0=\G_{\cH}(x')\geq \G_{\cH}(x^{f(H')})$, implying again by our assumption that $\T_{\cH}(x^{f(H')})=0$. Thus, $H'\subseteq supp(x)=S$ must intersect all hyperedges of $\cH_S$. Since this argument works for an arbitrary subset $S\subseteq V$ with $\cH_S\neq\emptyset$, property \eqref{intersection} follows.

For the other direction assume $\cH$ satisfies property \eqref{intersection}, and consider a position $x\in\ZZP^V$ for which $\T_{\cH}(x)>0$. Then $\cH_{supp(x)}\neq\emptyset$, and thus by property \eqref{intersection} we have a hyperedge $H\in\cH_{supp(x)}$ that intersects all other hyperedges of this induced subhypergraph, that is for which $\T_{\cH}(x^{f(H)})=0$, implying $\G_{\cH}(x^{f(H)})=0$ by \eqref{trivi-UB}. Since $x\to x^{f(H)}$ is an $H$-move, $\G_{\cH}(x)\neq 0$ is implied by the definition of the SG function. Since this follows for all positions $x$ with $\T_{\cH}(x)>0$, we can conclude that $\P_{\cH}\subseteq \Z_{\cH}$, as claimed.
\qed

\begin{corollary}\label{c111}
Condition \eqref{intersection} is necessary for a hypergraph to be Tetris.\qed
\end{corollary}

\smallskip

The following example demonstrates that condition \eqref{intersection} alone is not enough, generally, to guarantee that a hypergraph is Tetris, or equivalently by (ii) of Lemma \ref{l2-gap}, to ensure the equality of the SG and Tetris functions.

\begin{figure}

\begin{center}
\begin{tikzpicture}

\def \radius {3cm}
\def \margin {6} 

\foreach \s in {0,...,8}
{
  \node[draw,circle] (\s) at ({90 - 40 * \s}:\radius) {$\s$};


}

\draw[line width=1mm,dotted,blue] (1) -- (2) -- (3) -- (1);

\draw[line width=1mm,red] (0) -- (1) -- (4) -- (6) -- (0);

\end{tikzpicture}

\caption{A hypergraph $\cH$ on the ground set $V=\ZZ_9$, with hyperedges $T_i=\{i,i+1,i+2\}$ and $F_i=\{i,i+1,i+4,i+6\}$ for $i\in \ZZ_9$, where additions are modulo $9$, that is,
$\cH=\{T_i,~ F_i\mid i\in \ZZ_9\}$. The figure shows $T_1$ (dotted, blue) and $F_0$ (solid, red.) This hypergraph has an intersecting hyperedge, but is not Tetris. \label{fig-1}}

\end{center}

\end{figure}
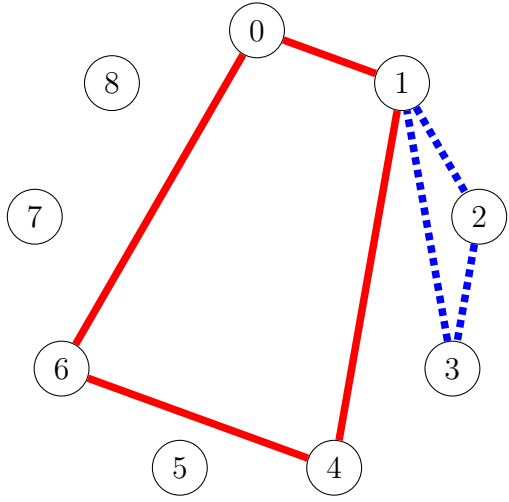

\begin{lemma}
The hypergraph defined in Figure \ref{fig-1} satisfies \eqref{intersection}, but does not have $\G_\cH=\T_\cH$.
\end{lemma}

\proof
To see this claim let us set $T_j$ and $F_j$ for $j\in\ZZ_9$ as in the caption of Figure \ref{fig-1}, where additions are modulo $9$.
Then let us observe first that
\[
T_j\cap F_i\neq\emptyset \text{ and } F_j\cap F_i\neq \emptyset \text{ for all } i,j\in\ZZ_9.
\]
An easy analysis show that $\cH$ satisfy condition \eqref{intersection}.

On the other hand, for the position $x=(1,1,...,1)\in\ZZ_+^V$ we have $\T_\cH(x)=3$.
Furthermore, $\T_\cH(x-\chi(T_j))=2$ and $\T_\cH(x-\chi(F_j))=0$ for all $j\in\ZZ_9$.
Thus, there exists no move $x\to x'$ with $\T_\cH(x')=1$, which by (iii) of Lemma \ref{l2-gap} implies that
$\G_\cH\neq \T_\cH$.
\qed

\bigskip

Note that the above example is a hypergraph with $\dim(\cH)=4$. We will show
later as claimed in Theorem \ref{t3} that for hypergraphs of dimension $3$ or less there are no such examples.

It is interesting to note that
for hypergraphs of dimension $2$ condition \eqref{intersection} can be substantially simplified.

\begin{lemma}\label{l000}
Assume that $\cH\subseteq 2^V$ is a hypergraph of $\dim(\cH)=2$ and
such that it has at least one edge $H\in\cH$ with $|H|=2$.
Then $\cH$ satisfies \eqref{intersection} if and only if there is a hyperedge $H\in\cH$
such that $H\cap H'\neq\emptyset$ for all $H'\in\cH$.
\end{lemma}

\proof
Figures below show the possible structure of such (hyper)graphs.
On the left $H$ is a singleton (red circle), while on te right it is a $2$-element set (red edge).
Circles in both pictures indicate possible singletons ($1$-element hyperedges.)
\qed

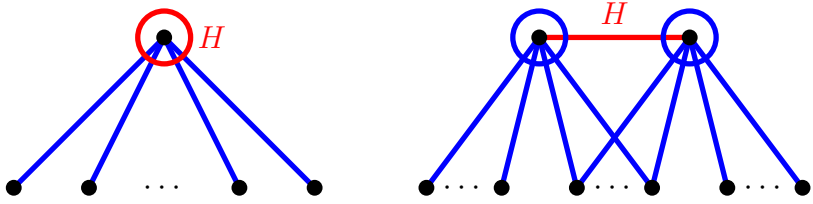
\begin{tikzpicture}[scale=0.5]

\coordinate (a) at (0,0);
\coordinate (b) at (2,0);
\coordinate (c) at (6,0);
\coordinate (d) at (8,0);
\coordinate (r) at (4,4);

\foreach \i in {a,b,c,d} \draw[line width=2pt, blue] (r) -- (\i);

\draw[red, line width=2pt] (4,4) circle (20pt);

\foreach \i in {a,b,c,d,r} \fill[black] (\i) circle (6pt);

\node[right of=r,node distance=18pt] {\color{red} $H$};

\node at (4,0) {$\cdots$};

\end{tikzpicture}
\hspace*{1cm}
\begin{tikzpicture}[scale=0.5]

\coordinate (a) at (0,0);
\coordinate (b) at (2,0);
\coordinate (c) at (4,0);
\coordinate (d) at (6,0);
\coordinate (e) at (8,0);
\coordinate (f) at (10,0);
\coordinate (u) at (3,4);
\coordinate (v) at (7,4);

\draw[line width=2pt, red] (u) -- (v) node[draw=none,fill=none,midway,above] {\color{red} $H$};

\foreach \i in {a,b,c,d} \draw[line width=2pt, blue] (u) -- (\i);

\foreach \i in {c,d,e,f} \draw[line width=2pt, blue] (v) -- (\i);

\draw[blue,line width=2pt] (3,4) circle (20pt);
\draw[blue,line width=2pt] (7,4) circle (20pt);

\foreach \i in {a,b,c,d,e,f,u,v} \fill[black] (\i) circle (6pt);

\node at (1,0) {$\cdots$};
\node at (5,0) {$\cdots$};
\node at (9,0) {$\cdots$};

\end{tikzpicture}

\bigskip

\subsection{A sufficient condition for $\T_\cH=\G_\cH$}\label{s-sufficient-dim}

Let us first recall some basic properties of Tetris functions.

\begin{lemma}\label{l1-interval}
Consider a hypergraph $\cH\subseteq 2^V$, a hyperedge $H\in\cH$, and a position $x\in\ZZP^V$ such that $x\geq \chi(H)$. Then for every integer value $\T_{\cH}(x^{f(H)})\leq v\leq \T_{\cH}(x^{s(H)})$ we have an $H$-move $x\to x'$ such that $\T_{\cH}(x')=v$.
\end{lemma}

\proof
We can decrease the components $x^{s(H)}_i$, $i\in H$ in an arbitrary order, subtracting one in each step, until we get $x^{f(H)}$. Every time the Tetris value can decrease by at most one. Since $x^{s(H)}_i<x_i$ for all $i\in H$, all the positions encountered in the above process can be reached from $x$ by a single $H$-move.
\qed

\smallskip

Given a hypergraph $\cH\subseteq 2^V$, and a position $x\in\ZZP^V$, let us call an integer vector $m\in\ZZP^\cH$ an $x$-vector, if
\begin{equation}\label{x-vector}
\sum_{H\in\cH} m_H \chi(H) ~\leq~ x ~~~\text{and}~~~ \sum_{H\in\cH} m_H ~=~ \T_{\cH}(x).
\end{equation}
Let us denote by $M(x)\subseteq \ZZP^\cH$ the family of $x$-vectors. Let us further define
\begin{equation}\label{x-pack}
\cH^{x-pack} ~=~ \{ H\in \cH \mid \exists m\in M(x) ~\text{s.t.}~ m_H>0\},
\end{equation}
that is $\cH^{x-pack}$ is the subfamily of $\cH$ of those hyperedges that participate with a positive multiplicity in some maximal $\T_{\cH}(x)$-packing of $\cH$. Every vector $m\in M(x)$ corresponds to such a maximal $\T_{\cH}(x)$-packing of $\cH$.

Let us consider $\cH_{supp(x)}$ the subhypergraph induced by the support of $x$, and define a subhypergraph of $\cH_{supp(x)}$ as
\begin{equation}\label{x-all}
\cH^{x-all} ~=~ \{H\in\cH_{supp(x)}\mid \forall H'\in \cH_{supp(x)}: H\cap H'\neq \emptyset\},
\end{equation}
consisting of those hyperedges that intersect all others in this subhypergraph.

\smallskip

\begin{lemma}\label{l3-basic}
Consider a hypergraph $\cH\subseteq 2^V$ that satisfies condition \eqref{intersection}, and a position $x\in\ZZP^V$ such that $\T_{\cH}(x)>0$. Then we have $\cH^{x-all}\neq \emptyset$ and $\cH^{x-pack} \neq \emptyset$. Furthermore, we have
\begin{itemize} \itemsep0em
\item[\rm{(i)}] $\T_{\cH}(x) ~>~ \T_{\cH}(x^{s(H)})~\geq~ \T_{\cH}(x^{f(H)}) ~\geq~ 0$ for all $H\in \cH_{supp(x)}$;
\item[\rm{(ii)}] $\T_{\cH}(x^{s(H)}) ~\geq~ \T_{\cH}(x)-|H|$ for all $H\in \cH_{supp(x)}$;
\item[\rm{(iii)}] $\T_{\cH}(x^{s(H)}) ~=~ \T_{\cH}(x)-1$ for all $H\in \cH^{x-pack}$;
\item[\rm{(iv)}] $\T_{\cH}(x^{f(H)}) ~=~ 0$ if and only if  $H\in \cH^{x-all}$;
\item[\rm{(v)}] $\T_{\cH}(x) ~\geq~ \T_{\cH}(x-\chi(\{k\}))~\geq~ \T_{\cH}(x)-1$ for all $k\in supp(x)$.
\end{itemize}
\end{lemma}

\proof
Trivial by the definitions.
\qed

\bigskip

We are now ready to prove Theorem \ref{t3}.

\subsubsection*{Proof of Theorem \ref{t3}.}

By Corollary \ref{c111}, condition \eqref{intersection} is necessary for a hypergraph of any dimension to be Tetris.
We prove next that for hypergraphs of dimension at most $3$ condition \eqref{intersection} is also sufficient.

The proof is indirect. Assume that there exists a hypergraph $\cH\subseteq 2^V$ of $\dim(\cH)\leq 3$ satisfying condition \eqref{intersection} such that $\G_\cH \neq \T_\cH$.
By Lemma \ref{l2-gap} this implies the existence of a position $x\in \ZZP^V$ and a value $\T_{\cH}(x)>v\geq 0$ such that there exists no move $x\to x'$ with $\T_{\cH}(x')=v$.
Since condition \eqref{intersection} applies to all induced subhypergraphs, we can assume without any loss of generality that
\begin{equation}\label{e10-supp}
V=supp(x).
\end{equation}
Then, by Lemma \ref{l1-interval} it follows that for all $H\in\cH$ we must have
\begin{subequations}
\begin{align}
&\text{either}\hspace*{-2cm} & \T_{\cH}(x^{s(H)}) &\leq v-1, \label{e-5a}\\
&\text{or}\hspace*{-2cm} & \T_{\cH}(x^{f(H)}) &\geq v+1. \label{e-5b}
\end{align}
\end{subequations}
By (i) of Lemma \ref{l3-basic} we cannot have both \eqref{e-5a} and \eqref{e-5b} hold for a hyperedge $H\in\cH$. Thus, the above defines a unique partition of the hyperedges of $\cH$:
\begin{align}
\begin{split}
\cH_1 &=~\{H\in\cH\mid \T_{\cH}(x^{s(H)}) \leq v-1\} ~~~\text{and}\\
\cH_2 &=~\{H\in\cH\mid \T_{\cH}(x^{f(H)}) \geq v+1\}.
\end{split}
\end{align}

Thus, for $H\in\cH_1$ we get by (ii) of Lemma \ref{l3-basic} that
\[
\T_{\cH}(x)-3 ~\leq~ \T_{\cH}(x^{s(H)}) ~\leq~ v-1,
\]
while for $H\in \cH_2$ we get by (i) and (iii) of Lemma \ref{l3-basic} that
\[
\T_{\cH}(x)-1 ~\geq~ \T_{\cH}(x^{s(H)}) ~\geq~ \T_{\cH}(x^{f(H)}) ~\geq~ v+1.
\]
These inequalities together imply that we must have $v=\T_{\cH}(x)-2>0$, and that
\begin{subequations}
\begin{align}
&\text{for}~ H\in\cH_1 ~\text{we have} & \T_{\cH}(x^{s(H)}) &= \T_{\cH}(x)-3, \text{ and}\label{e-either1}\\
&\text{for}~ H\in\cH_2 ~\text{we have} & \T_{\cH}(x^{f(H)}) &= \T_{\cH}(x)-1. \label{e-or1}
\end{align}
\end{subequations}

\bigskip

The next series of claims help us to prove that we must have $\T_{\cH}(x)=3$, and that we have $x_i=1$ for all $i\in H\in\cH_1$.

\bigskip

\begin{lemma}\label{l11}
We have $\cH^{x-all}\subseteq \cH_1$.
\end{lemma}

\proof
For all $H\in\cH^{x-all}$ we have by definition $\T_{\cH}(x^{f(H)})=0 < \T_{\cH}(x)-1$. Thus, $H\in\cH_1$ follows.
\qed

\begin{lemma}\label{l11b}
For all $H\in\cH_1$ we have $|H|=3$.
\end{lemma}

\proof
The claim follows by the definition of $\cH_1$, (ii) of Lemma \ref{l3-basic}, and the assumption that $\dim(\cH)\leq 3$.
\qed

\begin{lemma}\label{l12}
We have $\cH^{x-pack}= \cH_2$.
\end{lemma}

\proof
By definition, for all $H\in\cH^{x-pack}$ we have $\T_{\cH}(x^{s(H)})=\T_{\cH}(x)-1>\T_{\cH}(x)-3$, implying $H\in\cH_2$. For $H\in \cH_2$ by (i) of Lemma \ref{l3-basic} it follows that $\T_{\cH}(x)>\T_{\cH}(x^{s(H)})\geq \T_{\cH}(x^{f(H)}) = \T_{\cH}(x)-1$, implying $\T_{\cH}(x^{s(H)}) = \T_{\cH}(x)-1$. Let us choose an arbitrary $m\in M(x^{s(H)})$ and define $m'_H=m_H+1$ and $m'_{H'}=m_{H'}$ for all $H'\neq H$. Then we have $m'\in M(x)$ and $m'_H>0$ implying $H\in\cH^{x-pack}$ by \eqref{x-pack}.
\qed

\begin{lemma}\label{l13}
For all $m\in M(x)$ and $H\in\cH$ we have $m_H\leq 1$.
\end{lemma}

\proof
If $m_H\geq 2$ for some $H\in\cH$, then for position $x'=x-2\chi(H)$ we have that $\T_{\cH}(x')=\T_{\cH}(x)-2$ and $x\to x'$ is a move, contradicting our assumption that there exists no such move.
\qed

\begin{lemma}\label{l14}
For all $H_1\in\cH^{x-all}$ and $H_2\in\cH^{x-pack}(=\cH_2)$ we have $|H_1\cap H_2|=1$.
\end{lemma}

\proof
Let us assume indirectly that $|H_1\cap H_2|\geq 2$. By Lemma \ref{l11} we have that $|H_1|=3$. Assume w.l.o.g. that $H_1=\{i,j,k\}$ and  $\{i,j\}\subseteq H_2$. Let us then define  position $x'$ by $x'_\ell=x_\ell$ for $\ell\not\in\{i,j\}$ and $x'_\ell=x_\ell-1$ for $\ell\in\{i,j\}$. Then we have $x'\geq x^{s(H_2)}$, implying
\[
\T_{\cH}(x')\geq \T_{\cH}(x^{s(H_2)})=\T_{\cH}(x)-1
\]
by the monotonicity of $\T_{\cH}$, (iii) of Lemma \ref{l3-basic}, and Lemma \ref{l12}.
Furthermore, we have $x'-\chi(\{k\}) \leq x^{s(H_1)}$ implying by (v) of Lemma \ref{l3-basic} that
\[
\T_{\cH}(x')-1 ~\leq~ \T_{\cH}(x'-\chi(\{k\})) ~\leq~ \T_{\cH}(x^{s(H_1)}).
\]
From the above $\T_{\cH}(x^{s(H_1)})~\geq~ \T_{\cH}(x)-2$ follows, contradicting \eqref{e-either1}. This contradiction proves that we must have $|H_1\cap H_2|\leq 1$, while the definition of $\cH^{x-all}$ implies $H_1\cap H_2\neq\emptyset$, concluding the proof of our claim.
\qed

For a multiplicity vector $m\in M(x)$ let us associate the corresponding position $x(m)$ defined by
\begin{equation}
x(m) ~=~ \sum_{H\in\cH}m(H)\chi(H).
\end{equation}

\begin{lemma}\label{l15}
For all $m\in M(x)$ and $i\in H^*\in\cH^{x-all}$ we have $x(m)_i=x_i$.
\end{lemma}

\proof
Clearly, we must have $x(m)\leq x$ for all $m\in M(x)$, by the definition of $M(x)$.
Assume indirectly that there exists $m\in M(x)$ an index $i\in H^*=\{i,j,k\}$ such that $x(m)_i<x_i$.
Then we have $x(m)\leq x-\chi(\{i\})$, implying by (v) of Lemma \ref{l3-basic} that
\[
\T_{\cH}(x) ~\geq~ \T_{\cH}(x-\chi(\{i\})) ~\geq~ \sum_{H\in\cH}m(H) ~=~ \T_{\cH}(x),
\]
from which $\T_{\cH}(x-\chi(\{i\}))= \T_{\cH}(x)$ follows. Thus, again by (v) of Lemma \ref{l3-basic}, we would get
\[
\T_{\cH}(x^{s(H^*)}) ~=~ \T_{\cH}((x-\chi(\{i\}))-\chi(\{j,k\})) ~\geq~ \T_{\cH}(x-\chi(\{i\}))-2 ~=~ \T_{\cH}(x)-2,
\]
contradicting \eqref{e-either1} and Lemma \ref{l11}. This contradiction proves our claim.
\qed

\begin{corollary}\label{c15}
For all $H^*\in\cH^{x-all}$ we have $\T_{\cH}(x)= \sum_{i\in H^*}x_i$.
\end{corollary}

\proof
Applying Lemma \ref{l14} for an $m\in M(x)$, and noting that $m(H)>0$ implies $H\in\cH_2$ by Lemma \ref{l12}, we can write
\begin{align*}
\T_{\cH}(x) &=~ \sum_{H\in\cH}m(H)\\
&=~ \sum_{H\in\cH_2}m(H)\\
&=~ \sum_{H\in\cH_2}m(H)|H\cap H^*| \\
&=~ \sum_{i\in H^*}\sum_{H\in\cH_2\atop H\ni i}m(H)\\
&=~ \sum_{i\in H^*} x(m)_i\\
&=~ \sum_{i\in H^*} x_i,
\end{align*}
where the last equality follows by Lemma \ref{l15}.
\qed

\begin{lemma}\label{l16}
For all $H^*\in\cH^{x-all}$ and all $i\in H^*$ we have $x_i=1$.
\end{lemma}

\proof
Let us fix a hyperedge $H^*=\{i,j,k\}\in\cH^{x-all}$ and note that Lemmas \ref{l14} and \ref{l15} imply the existence of a hyperedge $H_2\in \cH_2$ with $H_2\cap H^*=\{i\}$. Let us then consider an arbitrary multiplicity vector $m\in M(x^{f(H_2)})$. Let us note that for all $H\in\cH$ with $m(H)>0$ we must have $H\subseteq supp(x^{f(H_2)})\subseteq supp(x)$, and thus $H\cap (H^*\setminus H_2)\neq\emptyset$ by the definition of $\cH^{x-all}$. Thus, using \eqref{e-or1} we can write
\begin{align*}
\T_{\cH}(x)-1~=~ \T_{\cH}(x^{f(H_2)}) &=~ \sum_{H\in\cH}m(H) \\
&\leq~ \sum_{H\in\cH}m(H)|H\cap (H^*\setminus H_2)|\\
&=~ x(m)_j+x(m)_k\\
&\leq~ x^{f(H_2)}_j+x^{f(H_2)}_k\\
&=~ x_j+x_k\\
&=~ \T_{\cH}(x) - x_i.
\end{align*}
From the above $x_i\leq 1$ follows, while $H^*\subseteq supp(x)$ implies $x_i\geq 1$.
\qed

\begin{corollary}\label{c16}
We have $\T_{\cH}(x)=3$.
\end{corollary}

\proof
Corollary \ref{c15} and Lemma \ref{l16} imply $\T_{\cH}(x)=3$.
\qed

\begin{corollary}\label{c16b}
We have $\cH^{x-all}=\cH_1$.
\end{corollary}

\proof
By Corollary \ref{c16} and \eqref{e-either1} we have $\T_{\cH}(x^{s(H)})=0$ for every $H\in\cH_1$, implying $\cH_1\subseteq \cH^{x-all}$. Thus the claim follows by Lemma \ref{l11}.
\qed

\bigskip

As a consequence of the above, we can restate  \eqref{e-either1} - \eqref{e-or1} as follows:
\begin{subequations}
\begin{align}
\T_{\cH}(x^{s(H)}) &= 0 & \forall H\in\cH_1=\cH^{x-all}, \label{e-either0}\\
\T_{\cH}(x^{f(H)}) &= 2 & \forall H\in\cH_2=\cH^{x-pack}. \label{e-or2}
\end{align}
\end{subequations}
Furthermore, by Lemma \ref{l16} and Corollary \ref{c16b} we have
\begin{equation}\label{e-degree}
x_i=1 ~~~\text{for all}~~~ i\in \bigcup_{H\in\cH_1}H.
\end{equation}

\begin{lemma}\label{l16.5}
For every $H\in\cH_1$ and for every $i\in H$ there exists $H'\in\cH_2$ such that $H\cap H'=\{i\}$.
\end{lemma}

\proof
Since $\T_{\cH}(x)=3$ by Corollary \ref{c16}, the equalities in \eqref{e-degree} and Lemma \ref{l12} imply the claim.
\qed

\bigskip

In the rest of the proof we show that $\cH_1$ and $\cH_2$ have some special structure, from which we can derive a contradiction at the end.
To this end we show first that $\cH_1$ includes three hyperedges such that any two of those intersect in exactly one point.

\begin{lemma}\label{l17}
For all $H^*\in\cH_1$ and $i\in H^*$ there exists $H^{**}\in\cH_1$ such that $i\not\in H^{**}$.
\end{lemma}
\proof
By Lemma \ref{l11b} we have $|H^*|=3$, and therefore we must have a point $j\in H^*\setminus\{i\}$.
Lemma \ref{l16.5} imply the existence of a hyperdege $H\in\cH_2$ such that $H^*\cap H=\{j\}$, and therefore $i\not\in H$. This implies $\cH_{V\setminus\{i\}}\neq\emptyset$, since this induced subhypergraph contains $H$. Therefore, by condition \eqref{intersection} there exists a hyperedge $H^{**}\in \cH_{V\setminus\{i\}}$ that intersects all others in this induced subhypergraph.
Consequently, for all hyperedges $H'\in\cH$ such that $H'\cap H^{**}=\emptyset$ we must have $i\in H'$. Therefore, $\T_{\cH}(x^{f(H^{**})})\leq x_i=1$. By the definition of $\cH_2$ and \eqref{e-or2} we get $H^{**}\in\cH_1$, as claimed.
\qed

\begin{lemma}\label{l18}
Consider $H_1,H_2\in\cH_1$ such that $i\in H_1\cap H_2$.
Then there exist no $H_3\in \cH_1$ such that $H_3\subseteq (H_1\cup H_2)\setminus\{i\}$.
\end{lemma}

\proof
By Lemma \ref{l16.5} there exists a hyperedge $H'\in\cH_2$ such that $H_1\cap H'=\{i\}$.
By Lemma \ref{l14} we also must have $|H_2\cap H'|=1$, thus by $i\in H_2$ we get $H'\cap ((H_1\cup H_2)\setminus \{i\})=\emptyset$.
Since by Lemma \ref{l14} we must have $H'\cap H_3\neq\emptyset$ for all $H_3\in\cH_1$, the claimed relation is implied.
\qed

\begin{lemma}\label{l19}
There exists hyperedges $H_1,H_2\in\cH_1$ such that $|H_1\cap H_2|=1$.
\end{lemma}

\proof
By Lemma \ref{l17} we have $|\cH_1|\geq 2$, and no point belongs to all edges of $\cH_1$.
Since $\cH_1=\cH^{x-all}$ by Corollary \ref{c16}, it is an intersecting family.
Since $\dim(\cH_1)=3$, any two distinct hyperedges of $\cH_1$ intersect in one or two points.

Assume indirectly that any two (distinct) hyperedges of $\cH_1$ intersect in two points.
Pick arbitrary two hyperedges of $\cH_1$, say $H_1=\{i,j,k\}$ and $H_2=\{i,j,\ell\}$, and let $H_3\in\cH_1$ such that $i\not\in H_3$.
By Lemma \ref{l17} such an $H_3$ exists.
Then, $|H_1\cap H_3|\geq 2$ and $|H_2\cap H_3|\geq 2$ together with $i\not\in H_3$ imply $H_3=\{j,k,\ell\}$,
that is that $H_3\subseteq (H_1\cup H_2)\setminus \{i\}$, contradicting Lemma \ref{l18}.
This contradiction proves the claim.
\qed

\begin{lemma}\label{l20}
There exists hyperedges $H_1,H_2,H_3\in\cH_1$ such that $|H_p\cap H_q|=1$ for all $1\leq p<q\leq 3$.
\end{lemma}

\proof
By Lemma \ref{l19} we have $H_1,H_2\in\cH_1$ such that $H_1\cap H_2=\{i\}$ for some $i\in V$.
Then by Lemma \ref{l17} there exists $H_3\in\cH_1$ such that $i\not\in H_3$.
We also have $H_3\not\subseteq (H_1\cup H_2)\setminus\{i\}$ by Lemma \ref{l18}.
Thus the claim follows.
\qed

\bigskip

\begin{corollary}\label{c20}
Thus, there exist six distinct points $X=\{a,b,c,d,e,f\}\subseteq V$ such that $H_1=\{a,b,f\}$, $H_2=\{b,c,d\}$ and $H_3=\{c,a,e\}$ are all hyperedges in $\cH_1$.\qed
\end{corollary}

\bigskip

We show next that $\cH_2$ has also a special form with respect to these six points.

\begin{lemma}\label{l21}
For all $H\in\cH_2$ we have one of the following: $\{a,d\}\subseteq H$, $\{b,e\}\subseteq H$, or $\{c,f\}\subseteq H$.
\end{lemma}

\proof
By Lemmas  \ref{l12}, \ref{l14}, and Corollary \ref{c16b} we have $|H\cap H_p|=1$ for all $p=1,2,3$. Then either $H$ has the form as claimed, or $H=\{d,e,f\}$. In the latter case however, let us consider $H'\in\cH_2$ such that $H'\cap H=\emptyset$. Such an $H'$ must exist by the facts $\T_{\cH}(x)=3$ and $H\in \cH_2=\cH^{x-pack}$. This set also must  intersect $H_p$, $p=1,2,3$ in exactly one point, however this is now impossible without intersecting $H$, too. Thus, only the claimed forms remain feasible for sets of $\cH_2$.
\qed

\begin{corollary}\label{c21}
Thus, using $\ga=\{a,d\}$, $\gb=\{b,e\}$ and $\gc=\{c,f\}$, we can conclude that the subhypergraphs
\begin{align*}
\cH_{2,\ga} &=\{H\in \cH_2\mid \ga\subseteq H\}, \\
\cH_{2,\gb} &=\{H\in \cH_2\mid \gb\subseteq H\}, \text{ and}\\
\cH_{2,\gc} &=\{H\in \cH_2\mid \gc\subseteq H\}
\end{align*}
form a partition of $\cH_2$. Furthermore, none of these families are empty.
\end{corollary}

\proof
The first claim follows directly by Lemma \ref{l21}. By \eqref{e-degree} we have $x_a=x_b=x_c=x_d=x_e=x_f=1$, and thus for any $m\in M(x)$ and $\mu\in\{\ga,\gb,\gc\}$ we must have $\sum_{H\in\cH_{2,\mu}} m(H)\leq 1$ by Lemma \ref{l15}. On the other hand we have $\T_{\cH}(x)=3$ by Corollary \ref{c16}, and thus for all $m\in M(x)$ and for all $\mu\in\{\ga,\gb,\gc\}$ we must have a hyperedge $H\in\cH_{2,\mu}$ with $m(H)=1$, completing the proof of the claim.
\qed

\begin{lemma}\label{l22}
These exists no hyperedge $H\in\cH_1$ that would contain $\mu$ for $\mu\in\{\ga,\gb,\gc\}$.
\end{lemma}

\proof
Assume indirectly that e.g., $H=\{a,d,u\}\in\cH_1$. Then by Lemma \ref{l14} we must have $u\in H'$ for all
$H'\in\cH_{2,\gb}\cup\cH_{2,\gc}$, and thus, in particular, $u\not\in X$. Since $\T_{\cH}(x)=3$, we must have $x_u\geq 3$. Let us then consider the $H$-move $x\to x'$, where $x'_i=x_i$ for $i\not\in H$, $x'_a=0$, $x'_d=0$, and $x'_u=1$. Then all hyperedges of  $\cH_2$ that are subsets of $supp(x')$ contain $u$, and thus we must have $\T_{\cH}(x')=1$, contradicting \eqref{e-either0}.
\qed

Let us next introduce $N_\mu=\bigcup\{H\setminus \mu\mid H\in\cH_{2,\mu}\}$ for $\mu\in\{\ga,\gb,\gc\}$. Note that these sets are disjoint from $X=H_1\cup H_2\cup H_3$, defined in Corollary \ref{c20}, by Lemma \ref{l14}.

\begin{lemma}\label{l23}
Let $\mu,\nu\in\{\ga,\gb,\gc\}$, $\mu\neq \nu$ and consider two sets $H\in\cH_{2,\mu}$ and $H'\in\cH_{2,\nu}$ that are disjoint $H\cap H'=\emptyset$. Then there exists a hyperedge $H''\in\cH_{2,\mu}\cup \cH_{2,\nu}$ that intersects both $H$ and $H'$.
\end{lemma}

\proof
By condition \eqref{intersection} we must have a set $H''\subseteq H\cup H'$, $H''\in\cH$ that intersects all sets in the non-empty induced subhypergraph $\cH_{H\cup H'}$. If $H''\in\cH_1$, then $|H''|=3$ by Lemma \ref{l11b}, and thus we must have either $|H''\cap H|\geq 2$ or $|H''\cap H'|\geq 2$, contradicting Lemma \ref{l14}. Thus we must have $H''\in\cH_2$, and therefore $H''\in \cH_{2,\mu}\cup \cH_{2,\nu}$ by Lemma \ref{l21}, as claimed.
\qed

\begin{corollary}\label{c23}
For $\mu,\nu\in\{\ga,\gb,\gc\}$, $\mu\neq \nu$, we either have $N_\mu\subseteq N_\nu$ or $N_\nu\subseteq N_\mu$.
\end{corollary}

\proof
If there are points $u\in N_\mu\setminus N_\nu$ and $v\in N_\nu\setminus N_\mu$, then by Lemma \ref{l23} we have either $\mu\cup\{v\}\in\cH_{2,\mu}$, or $\nu\cup\{v\}\in\cH_{2,\nu}$ contradicting $u\not\in N_\nu$ or $v\not\in N_\mu$.
\qed

\begin{lemma}\label{l24}
Let $\mu,\nu\in\{\ga,\gb,\gc\}$, $\mu\neq \nu$. Then, there exists no two distinct points $u,v\in V\setminus X$ such that all four sets $\mu\cup\{u\}$, $\mu\cup\{v\}$, $\nu\cup\{u\}$, and $\nu\cup\{v\}$ are hyperedges of $\cH$.
\end{lemma}

\proof
Assume indirectly that such points do exist.
Then by Lemma \ref{l22} these sets are all from $\cH_2$.
By condition \eqref{intersection} we must have a hyperedge $H\subseteq \{\mu\cup\nu\cup\{u,v\}$ in $\cH$ that intersects all these sets.
Since $H$ must intersect some of these four sets in two points, $H\in \cH_2$ holds by Lemma \ref{l14}. Then, by Corollary \ref{c21} we have $H\in \cH_{2,\mu}\cup\cH_{2,\nu}$. This is however impossible, since there exists no such subset of size at most $3$ that would either contain $\mu$ or $\nu$ and intersect all these fours sets.
\qed

\begin{corollary}\label{c24}
For all $\mu,\nu\in\{\ga,\gb,\gc\}$, $\mu\neq \nu$ we have $|N_\mu\cap N_\nu|\leq 1$.
\end{corollary}

\proof
Immediate from Lemma \ref{l24}.
\qed

\begin{corollary}\label{c25}
Up to a relabeling of the vertices, we have $N_\ga\subseteq N_\gb\subseteq N_\gc$, and $|N_\ga|\leq |N_\gb|\leq 1$.
\end{corollary}

\proof
Immediate by Corollaries \ref{c23} and \ref{c24}.
\qed

\begin{lemma}\label{l26}
At most one of $\ga$, $\gb$ and $\gc$ is a hyperedge of $\cH$.
\end{lemma}

\proof
If e.g., $\ga,\gb\in\cH$, then by property \eqref{intersection} we must have a hyperedge $H\in\cH$ such that $H\subseteq \ga\cup\gb$ and it intersects both $\ga$ and $\gb$. Since $\ga,\gb\in\cH_2$, by Lemmas \ref{l14} and \ref{l11b} we have that $H\not\in\cH_1$ and hence $H\in\cH_2$. Then, by Corollary \ref{c21} we must have $H\in\cH_{2,\ga}$ or $H\in\cH_{2,\gb}$. Since $H$ must intersect both $\ga$ and $\gb$, $|H|=3$ follows, from which we derive a contradiction by Lemma \ref{l14},
due to the structure of $\cH_1$ sets within the set $X$.
\qed

\begin{lemma}\label{l25}
$N_\ga\neq\emptyset$.
\end{lemma}

\proof
Assume indirectly that $N_\ga=\emptyset$. This implies that $\cH_{2,\ga}=\{\ga\}$. Let us now consider an arbitrary $m\in M(x)$. Since $\T_{\cH}(x)=3$ by Corollary \ref{c16}, we
must have hyperedges $H_\mu\in\cH_{2,\mu}$ for all $\mu\in\{\ga,\gb,\gc\}$ with $m(H_\mu)=1$ by \eqref{e-degree}. In particular, we must have $m(\ga)=1$ and $m(H)=1$ for some $H\in\cH_{2,\gb}$. Since $\ga\cap H=\emptyset$, by property \eqref{intersection} we must have a hyperedge $H'\in\cH$ that intersects both $\ga$ and $H$ such that $H'\subseteq \ga\cup H$. If $H'\in\cH_1$ then we get a contradiction by Lemma \ref{l14}. Thus we must have $H'\in \cH_2$. Then by Corollary \ref{c21} and the fact that $\cH_{2,\ga}=\{\ga\}$ we must have $H'\in\cH_{2,\gb}$. This contradicts the fact that $\ga$ is disjoint from all sets of $\cH_{2,\gb}$.
\qed

Thus by Corollary \ref{c25} and Lemma \ref{l25} we have $|N_\ga|=|N_\gb|=1$, that is for some $u\in V$ we have
$N_\ga=N_\gb=\{u\}\subseteq N_\gc$. Therefore we have $H=\gc\cup\{ u\}\in\cH_{2,\gc}$. Let $x'=x^{f(H)}$, and consider $m\in M(x')$. By Lemma \ref{l14} we have $m(H^*)=0$ for all $H^*\in\cH_1$. Furthermore, for any $H'\in\cH_{2}$ such that $u\in H'$ we also must have $m(H')=0$.
Consequently, only $H'\in\cH_{2,\ga}\cup\cH_{2,\gb}$, $u\not\in H'$ can have $m(H')=1$ (and not more by \eqref{e-degree}.)
Since by Lemma \ref{l26} at most one of $\ga$ and $\gb$ can belong to $\cH_2$, we must have $\T_{\cH}(x')\leq 1$ contradicting \eqref{e-or2}. This contradiction completes the proof of the theorem.
\qed

\bigskip

\subsection{Another sufficient condition for $\T_\cH=\G_\cH$}\label{s-sufficient-int}

We can strengthen condition \eqref{intersection} by requiring that any two hyperedges intersect. In other words, for all $H,H'\in\cH$ we have $H\cap H'\neq\emptyset$.

\begin{theorem}
If $\cH\subseteq 2^V$ is an intersecting hypergraph, then we have $\G_{\cH}(x)=\T_{\cH}(x)$ for all positions $x\in\ZZP^V$.
\end{theorem}

\proof
Let us consider an arbitrary position $x\in\ZZP^V$. If $\T_{\cH}(x)=0$, then the claim holds by definition. Assume that $\T_{\cH}(x)>0$ and consider a hyperedge $H\in\cH$ such that $\T_{\cH}(x-\chi_H)=\T_{\cH}(x)-1$. Such a hyperedge exists since $\T_{\cH}(x)>0$. Let us then consider the positions $x^{s(H)}$, and $x^{f(H)}$. By our choice of $H$ we have $\T_{\cH}(x^{s(H)})=\T_{\cH}(x)-1$. Since the hypergraph is intersecting, we also have $\T_{\cH}(x^{f(H)})=0$. Thus, by Lemma \ref{l1-interval} for all values $0\leq v\leq \T_{\cH}(x)-1$ there exists an $H$-move $x\to x'$ such that $\T_{\cH}(x')=v$. Since this holds for all positions, we get $\G_{\cH}(x)=\T_{\cH}(x)$ by Lemma \ref{l2-gap}. \qed

\bigskip

\subsection{Computing the Tetris Function}

\begin{theorem}\label{t-complexity}
Given a hypergraph $\cH\subseteq 2^V$ and a position $x\in\ZZP^V$,
computing $\T_\cH(x)$ is
\begin{itemize} \itemsep0em
\item[\rm{(i)}] NP-hard for intersecting hypergraphs;
\item[\rm{(ii)}] NP-hard for hypergraphs of dimension at most $3$;
\item[\rm{(iii)}] polynomial for hypergraphs of dimension at most $2$ (i.e., for graphs).
\end{itemize}
\end{theorem}

\proof
Let us consider an arbitrary hypergraph $\cH\subseteq 2^V$. Its matching number $\mu(\cH)$ is the maximum number of pairwise disjoint hyperedges of $\cH$, and is known to be NP-hard to compute. Let us then consider $w\not\in V$ and define $\cH^*=\{H\cup\{w\}\mid H\in\cH\}$. Furthermore, let us consider the position $x\in\ZZP^{V\cup\{w\}}$ defined by $x_i=1$ for $i\in V$ and $x_w=|\cH|$. Then $\cH^*$ is an intersecting hypergraph and we have $\T_{\cH^*}(x)=\mu(\cH)$.

If $\cH$ is of dimension $3$, and $x_i=1$ for all $i\in V$, then we have again $\T_\cH(x)=\mu(\cH)$.

Finally, if $\cH$ is of dimension at most $2$, then $T_\cH(b)$ for a position $b\in\ZZP^V$ is the so called $b$-matching number of the underlying graph and is known to be computable in polynomial time (see \cite{Edm65,Tutte}).
\qed

\begin{corollary}\label{c222}
Given a hypergraph $\cH\subseteq 2^V$ and a position $x\in\ZZP^V$,
computing $\G_\cH(x)$ is NP-hard, already for intersecting hypergraphs.
\end{corollary}

\proof
Since intersecting hypergraphs satisfy condition \eqref{intersection},
Theorem \ref{t3} implies $\T_\cH=\G_\cH$. Thus the claim follows by (i) of Theorem  \ref{t-complexity}.
\qed


\begin{thebibliography}{99}


\bibitem{Alb07}
M.H. Albert, R.J. Nowakowski, D. Wolfe,
Lessons in play: An introduction to combinatorial game theory,
second ed., A K Peters Ltd., Wellesley, MA, 2007.

\bibitem{BCG01-04}
E.R. Berlekamp, J.H. Conway, and R.K. Guy,
Winning ways for your mathematical plays,
vol.1-4, second edition, A.K. Peters, Natick, MA, 2001 - 2004.

\bibitem{BGHM15}
Endre Boros, Vladimir Gurvich, Nhan Bao Ho, and Kazuhisa Makino;
Extended complementary {\sc Nim},
RUTCOR Research Report, RRR-1-2015, Rutgers University;
http://arxiv.org/abs/1504.06926 . 

\bibitem{BGHMM15}
Endre Boros, Vladimir Gurvich, Nhan Bao Ho, Kazuhisa Makino, and Peter Mursic.
On the Sprague-Grundy function of exact $k$-Nim,
RUTCOR Research Report, RRR-2-2015, Rutgers University;
http://arxiv.org/abs/1508.04484 . 

\bibitem{BGHMM16}
Endre Boros, Vladimir Gurvich, Nhan Bao Ho, Kazuhisa Makino, and Peter Mursic,
On the Sprague-Grundy function of hypergraph combinations of impartial games,
manuscript.


\bibitem{Bou901}
C.L. Bouton, Nim, a game with a complete mathematical theory,
Ann. of Math., 2-nd Ser. 3 (1901-1902) 35--39.



\bibitem{Edm65}
J. Edmonds,
Paths, trees and flowers,
Canadian Journal of Mathematics, 17 (1965) 449–-467.

\bibitem{Gru39}
P.M. Grundy, Mathematics of games, Eureka 2 (1939) 6--8.

\bibitem{GS56}
P.M. Grundy and C.A.B. Smith,
Disjunctive games with the last player loosing,
Proc. Cambridge Philos. Soc., 52 (1956) 527--523.

\bibitem{JM80}
T.A. Jenkyns and J.P. Mayberry,
Int. J. of Game Theory 9 (1) (1980) 51--63,
The skeletion of an impartial game
and the Nim-Function of Moore's Nim$_k$.

\bibitem{Moo910}
E. H. Moore,
A generalization of the game called {\sc Nim},
Annals of Math., Second Series, 11:3 (1910) 93--94.

\bibitem{NM44}
J. von Neumann and O. Morgenstern,
Theory of games and economic behavior, Princeton University Press, 1944.



\bibitem{Smi66}
C. A. B. Smith, Graphs and composite games,
J. of Combinatorial theory 1  (1966) 51--81.

\bibitem{Spr35}
R. Sprague, \"Uber mathematische Kampfspiele,
Tohoku Math. J. 41 (1935-36) 438--444.

\bibitem{Spr37}
R. Sprague, \"Uber zwei abarten von nim, Tohoku Math. J. 43
(1937) 351--354.

\bibitem{Tutte}
W. Tutte, A short proof of the factor theorem for finite graphs,
Canadian Journal of Mathematics, 6 (1954) 347–-352.

\end{thebibliography}
\end{document}